\documentclass{article}
\usepackage{spconf,amsmath,graphicx,cite}

\title{Optimal ROC Curves from Score Variable Threshold Tests}
\name{Catherine Medlock (IEEE Student Member) and Alan Oppenheim (IEEE Life Fellow) \thanks{*This work was supported through the generosity of Analog Devices Inc., Texas Instruments Inc., and Bose Corporation.}}
\address{Massachusetts Institute of Technology \\
Department of Electrical Engineering and Computer Science \\
77 Massachusetts Ave., Cambridge MA, 02139}
\begin{document}
%
\maketitle
\begin{abstract}
The Receiver Operating Characteristic (ROC) is a well-established representation of the tradeoff between detection and false alarm probabilities in binary hypothesis testing. In many practical contexts ROC's are generated by thresholding a measured score variable -- applying score variable threshold tests (SVT's). In many cases the resulting curve is different from the likelihood ratio test (LRT) ROC and is therefore not Neyman-Pearson optimal. While it is well-understood that concavity is a {\bf necessary} condition for an ROC to be Neyman-Pearson optimal, this paper establishes that it is also a {\bf sufficient} condition in the case where the ROC was generated using SVT's. It further defines a constructive procedure by which the LRT ROC can be generated from a non-concave SVT ROC, without requiring explicit knowledge of the conditional PDF's of the score variable. If the conditional PDF's are known, the procedure implicitly provides a way of redesigning the test so that it is equivalent to an LRT.
\end{abstract}
\section{Introduction}
\label{sec:intro}

Receiver Operating Characteristic (ROC) curves have played a key role in the signal detection communities for many decades \cite{van1968detection,helstrom1968}. In those contexts, a family of ROC's is typically generated using threshold tests on the likelihood ratio (LRT's) of a measured value often called the score variable. Each curve displays the probability of detection (true positive rate) vs. probability of false alarm (false positive rate) for a fixed SNR and is parameterized by the LRT threshold value. It is well-known that LRT's are optimal in the Neyman-Pearson sense, i.e., for any upper bound on the probability of false alarm, the probability of detection is maximized. In a typical radar signal detection scenario, the score variable often used is the sampled signal level at the output of a matched filter. And in this scenario the conditional PDF's of the score variable are typically based on a mathematical model for the noise with their analytical forms used to perform LRT's.

It is widely known that an ROC curve generated using LRT's on any pair of conditional PDF's, assuming it is continuous, must be strictly concave, i.e., strict concavity is a {\bf necessary} condition for the Neyman-Pearson optimality of a continuous ROC curve. It is also known that if the likelihood ratio is an invertible function of the score variable for a given set of conditional PDF's, then threshold tests on the score variable (SVT's) and LRT's will result in identical ROC curves. Equivalently, in this case, the score variable is a sufficient statistic for the likelihood ratio.
 
The use of ROC curves has become increasingly prevalent in a very broad set of application areas including biostatistics and machine learning \cite{fawcett2006introduction,obuchowski2003receiver,zweig1993receiver,spackman1989signal,beck1986use,sainath2015convolutional,stein2005relationship}. In contrast to the problem of signal detection, in these contexts the score variable is typically a finely-tuned combination of many measurements and is therefore often less amenable to mathematical analysis and modeling. In principle, this does not preclude the use of LRT's to generate ROC curves, since in the absence of specific models of the conditional PDF's of the score variable, the PDF's can be estimated from histograms derived from training data. However, reliable estimation of probability densities from empirical data is well-known to be a difficult problem \cite{rosenblatt1956remarks,parzen1962estimation}. Estimation of the likelihood ratio from empirical data is even more so because small errors in the estimate of the denominator of the ratio can lead to large errors in the estimate of the ratio itself. SVT's are simpler to perform and less sensitive to poor estimates of the probability densities. Consequently SVT's are often used instead of LRT's.

This paper identifies and resolves the issue of when an ROC curve generated with SVT's is equivalent to the Neyman-Pearson optimal curve. Section \ref{sec:background} defines the terminology and notation used throughout the rest of the paper. Throughout the paper, the word ``optimal'' is used to mean ``Neyman-Pearson optimal.'' Section \ref{sec:convex} establishes that in the case where the ROC curve is generated using SVT's, concavity is not only a necessary condition but also a {\bf sufficient} condition for optimality. Section \ref{sec:convex} further establishes a constructive procedure for obtaining the optimal ROC curve of a given score variable from its SVT ROC curve, without requiring knowledge of the underlying conditional PDF's. If the conditional PDF's are known, or equivalently the SVT threshold value associated with each point on the SVT ROC curve is known, then the procedure also implicitly defines a way of redesigning the test so that it is optimal.

\section{Background} \label{sec:background}

\noindent We refer to an arbitrary score variable $S$ whose conditional PDF's under the hypotheses $H_0$ (the null hypothesis) and $H_1$ (the positive hypothesis) are $f_0(\cdot)$ and $f_1(\cdot)$, respectively, using the notation
	\begin{equation}
	S \sim f_0(\cdot), f_1(\cdot).
	\end{equation}
We assume for simplicity that all score variables are continuous, have continuous support, and have smooth LRT ROC curves.\footnote{This translates into the properties that the ratio $f_1(s)/f_0(s)$ takes on all values in some interval $[0,M)$ for $M>0$ and that the set ${\{s : f_1(s)/f_0(s) = \eta\}}$ is zero measure for all $\eta \geq 0$. These properties also ensure that if the SVT ROC curve is concave, then it must be strictly concave.} Score variables that do not satisfy these assumptions can be treated using similar methods to those in this paper, but are not included for the sake of brevity.

Denote by $\mathcal{D}$ the set of values of $S$ for which $H_1$ is chosen by a given decision rule. Each point on an ROC curve has coordinates
\begin{subequations}
\begin{alignat}{1}
P_F &= \int_{\mathcal{D}} ds \, f_0(s) \\[10pt]
P_D &= \int_{\mathcal{D}} ds \, f_1(s).
\end{alignat}
\end{subequations}
If the decision rule is an SVT with threshold $\gamma$, then $\mathcal{D}$ takes the form
\begin{equation}
\mathcal{D}_{\text{SVT}}(\gamma) = \{ s : s \geq \gamma \}.
\end{equation}
If it is an LRT with threshold $\eta \geq 0$, then $\mathcal{D}$ takes the form
\begin{equation} \label{eq:lrt_decision_region}
\mathcal{D}_{\text{LRT}}(\eta) = \{ s : f_1(s) / f_0(s) \geq \eta \}.
\end{equation}
Note that while $\mathcal{D}_{\text{SVT}}(\gamma)$ is always a contiguous region of the real line, in general $\mathcal{D}_{\text{LRT}}(\eta)$ consists of a set of disjoint intervals of $s$. Obviously, this means that for an arbitrary pair of conditional PDF's, the SVT and LRT ROC curves can be different.

\section{Optimality of a Concave SVT ROC Curve} \label{sec:concave}

\noindent A principal result of this paper is that if an ROC generated using SVT's on a score variable $S \sim f_0(\cdot), f_1(\cdot)$ is concave, then it is {\bf guaranteed} to be the optimal ROC curve for that score variable. To show that this is true, note that the SVT ROC curve can be written as
\begin{subequations} \label{fig:svtroc_0}
\begin{alignat}{1}
P_F^{\text{SVT}} &= \int_{\mathcal{D}_{\text{SVT}}(\gamma)} ds \, f_0(s) = \int_\gamma^\infty ds \, f_0(s) \\[10pt]
P_D^{\text{SVT}} &= \int_{\mathcal{D}_{\text{SVT}}(\gamma)} ds \, f_1(s) = \int_\gamma^\infty ds \, f_1(s).
\end{alignat}
\end{subequations}
where $-\infty < \gamma < +\infty$. Equation \eqref{fig:svtroc_0} can be simplified using the cumulative distribution functions (CDF's) $F_0(s) = \int_{-\infty}^s du \, f_0(u)$ and $F_1(s) = \int_{-\infty}^s du \, f_1(u)$ to obtain
\begin{subequations} \label{eq:svtroc_1}
\begin{alignat}{1}
P_F^{\text{SVT}} &= \phi_F(\gamma) = 1 - F_0(\gamma) \\
P_D^{\text{SVT}} &= \phi_D(\gamma) = 1 - F_1(\gamma)
\end{alignat}
\end{subequations}
where $-\infty < \gamma < +\infty$ and the functions $\phi_F(\cdot)$ and $\phi_D(\cdot)$ have been introduced for convenience later on. At the point on the curve corresponding to the SVT threshold $\gamma$, the derivative of the curve is
\begin{equation} \label{eq:derivative}
\left.\frac{d P_D^{\text{SVT}}}{d P_F^{\text{SVT}}}\right|_{P_F^{\text{SVT}} = \phi_F(\gamma)} = \frac{\phi_D^\prime(\gamma)}{\phi_F^\prime(\gamma)} = \frac{f_1(\gamma)}{f_0(\gamma)}.
\end{equation}
It follows that if the derivative of an SVT ROC curve is an invertible function of $P_F$, as is the case for a concave curve under the current assumptions, and if $\phi_F(\cdot)$ is strictly monotonic, as is the case for all score variables under the current assumptions, then the likelihood ratio is an invertible function of the score variable. This implies that the SVT ROC curve of the score variable is identical to its LRT ROC curve, proving the result.

It is important to point out that this result is different from the statement in \cite{torgersen1991comparison}, which says that given any concave curve with endpoints at $(0,0)$ and $(1,1)$, one can always construct a pair of conditional PDF's for which that curve is the optimal ROC curve. In that context, the curve and the PDF's are strictly abstract and need not have any connection to an actual binary hypothesis testing problem. On the other hand, the result stemming from Equation \eqref{eq:derivative} says that if the given curve is an ROC curve that was generated using SVT's on a {\bf specific} pair of PDF's and is strictly concave, then the curve is optimal {\bf for those PDF's}.

\section{Generation of the Optimal ROC Curve from a Non-Concave SVT ROC Curve} \label{sec:convex}

\noindent This section defines a procedure for constructing the optimal ROC curve of a score variable $S \sim f_0(\cdot), f_1(\cdot)$ directly from its SVT ROC curve, without requiring knowledge of the PDF's $f_0(\cdot)$ and $f_1(\cdot)$. It is assumed that the SVT ROC curve is not concave, since otherwise it would already be optimal according to Section \ref{sec:concave}. Mathematically, the objective is to be able to compute the functions
\begin{subequations} \label{eq:lrtroc_0}
\begin{alignat}{1}
P_F^{\text{LRT}} &= \int_{\mathcal{D}_{\text{LRT}}(\eta)} ds \, f_0(s) \\
P_D^{\text{LRT}} &= \int_{\mathcal{D}_{\text{LRT}}(\eta)} ds \, f_1(s)
\end{alignat}
\end{subequations}
directly from the SVT ROC curve of $S$ for any $\eta \geq 0$, assuming that $f_0(\cdot)$ and $f_1(\cdot)$ are unknown. This is equivalent to assuming that the SVT ROC curve as a whole is known, but the SVT threshold value associated with each individual point is unknown. Note that under these assumptions, the derivative $dP_D^{\text{SVT}}/dP_F^{\text{SVT}}$ can be computed as a function of $P_F^{\text{SVT}}$ but not as a function of the score variable $s$.

We will show that the following steps accomplish the objective.

\begin{figure}
	\centering
	\includegraphics[width=3in]{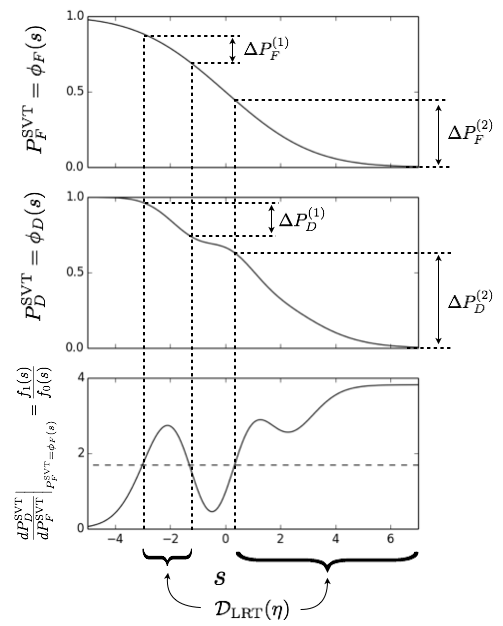}
	\caption{Procedure for generating optimal LRT ROC curve from SVT ROC curve.}
	\label{fig:repair}
\end{figure}

\begin{enumerate}
	\item Choose a value of $\eta \geq 0$.
	\item Identify the segments of the SVT ROC curve where the slope $dP_D^{\text{SVT}}/dP_F^{\text{SVT}}$ is greater than or equal to $\eta$. Record the changes in $P_F^{\text{SVT}}$ and the changes in $P_D^{\text{SVT}}$ associated with the endpoints of each segment. Denote these changes by $\Delta P_F^{(i)}$ and $\Delta P_D^{(i)}$ where $i$ is an index over segments.
	\item Compute $P_F^{\text{LRT}}$ and $P_D^{\text{LRT}}$ of the optimal curve by summing the $\Delta P_F^{(i)}$ and $\Delta P_D^{(i)}$ values,
	\begin{subequations}
		\begin{alignat}{1}
			P_F^{\text{LRT}} &= \sum_i \Delta P_F^{(i)} \\
			P_D^{\text{LRT}} &= \sum_i \Delta P_D^{(i)}.
		\end{alignat}
	\end{subequations}
	Equivalently, add the segments together end-to-end to compute the location of the point on the optimal curve.
\end{enumerate}
The underlying logic can be understood through Figure \ref{fig:repair}. The top and bottom graphs show the probability of false alarm $P_F^{\text{SVT}} = \phi_F(s)$ and the derivative $\left. dP_D^{\text{SVT}}/dP_F^{\text{SVT}}\right|_{{P_F}^{\text{SVT}} = \phi_F(s)}$ as functions of $s$, respectively. From Equation \ref{eq:derivative}, the bottom graph is equivalent to the likelihood ratio $f_1(s)/f_0(s)$. Note that these graphs are caricatures used only for visualization, since the procedure does not require explicit knowledge of either $P_F^{\text{SVT}}$ or $dP_D^{\text{SVT}}/dP_F^{\text{SVT}}$ as functions of $s$. It does require knowledge of $dP_D^{\text{SVT}}/dP_F^{\text{SVT}}$ as a function of $P_F^{\text{SVT}}$.

\begin{figure}
	\centering
	\includegraphics[width=2.5in]{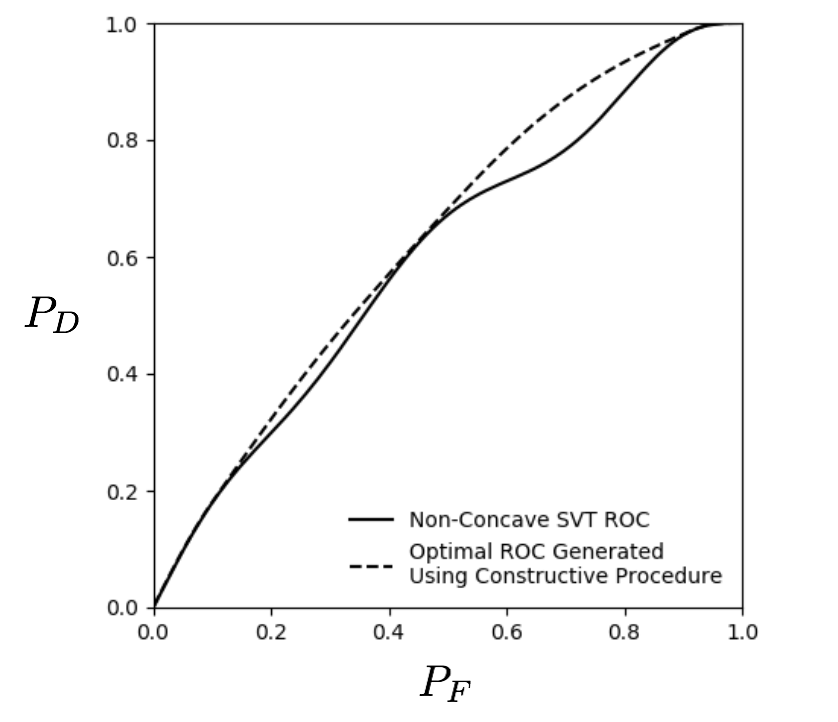}
	\caption{Non-concave SVT ROC and optimal ROC generated using the procedure given in the text.}
	\label{fig:roc_final}
\end{figure}

In the bottom graph, a fixed LRT threshold value $\eta \geq 0$ identifies multiple disjoint regions of $s$ for which $f_1(s)/f_0(s) \geq \eta$, and together these regions comprise $\mathcal{D}_{\text{LRT}}(\eta)$. Each individual region $i$ has the form $[a_i,b_i]$ for $a_i < b_i$ and corresponds to the segment in the SVT ROC curve with endpoints $(P_F^{\text{SVT}}, P_D^{\text{SVT}}) = (\phi_F(b_i),\phi_D(b_i))$ and $(P_F^{\text{SVT}}, P_D^{\text{SVT}}) = (\phi_F(a_i),\phi_D(a_i))$. The integrals of $f_0(\cdot)$ and $f_1(\cdot)$ over the region are then simply
\begin{subequations}
\begin{alignat}{1}
\int_{a_i}^{b_i} ds \, f_0(s) &= (1-F_0(a_i)) - (1-F_0(b_i)) \nonumber\\
&= \phi_F(a_i) - \phi_F(b_i) \\
\int_{a_i}^{b_i} ds \, f_1(s) &= (1-F_1(a_i)) - (1-F_1(b_i)) \nonumber\\
&= \phi_D(a_i) - \phi_D(b_i)
\end{alignat}
\end{subequations}
which are exactly the changes in $P_F^{\text{SVT}}$ and $P_D^{\text{SVT}}$ between the endpoints of the segment. Summing these changes over all regions corresponds to summing the integrals of $f_0(\cdot)$ and $f_1(\cdot)$ over each disjoint portion of $\mathcal{D}_{\text{LRT}}(\eta)$. The resulting optimal ROC curve made by varying $\eta$ over its entire range is shown in Figure \ref{fig:roc_final}.

Note that this procedure is different than the use of randomization to replace a convex region on an ROC curve by the straight line connecting its endpoints \cite{van1968detection,provost1997analysis}. In that case, a biased coin is flipped and the result dictates whether the decision region of the first endpoint is used or that of the second. It is straightforward to show that the effective probabilities of false alarm and detection then lie on the straight line between the endpoints. However, the resulting curve is not optimal. One way of seeing this is to observe that the LRT ROC curve of a continuous score variable, even in the absence of the assumptions made in this paper, can never have any linear regions -- it must either be continuous and strictly concave or discontinuous and strictly concave over each of its disjoint regions.

Consider now the case where the PDF's $f_0(\cdot)$ and $f_1(\cdot)$ are known. Equivalently, the SVT threshold $\gamma$ associated with each point on the SVT ROC curve is known. It may be desirable to redesign the test on the score variable to be an LRT instead of an SVT. For an LRT with threshold $\eta \geq 0$, this requires explicit knowledge of the regions of $s$ for which $f_1(s)/f_0(s) \geq \eta$, i.e., the regions of $s$ that comprise $\mathcal{D}_{\text{LRT}}(\eta)$. One way of obtaining the regions is to solve analytically for the values of $s$ for which the inequality holds. An alternative method is provided implicitly by the procedure given above and can be seen in the top graph in Figure \ref{fig:repair} -- Compute the SVT ROC curve of $f_0(\cdot)$ and $f_1(\cdot)$, plot the derivative of the curve as a function of the SVT threshold, then read the regions directly off the graph by checking where the derivative is greater than or equal to $\eta$. The redesigned test then decides $H_1$ for any value of $s$ that lies within $\mathcal{D}_{\text{LRT}}(\eta)$.

\section{Summary}

In summary, while it is well-understood that concavity is a {\bf necessary} condition for an ROC curve to be optimal for the score variable used to generate it, in this paper we have established that it is also a {\bf sufficient} condition in the case where the ROC curve was generated using SVT's. We then defined a procedure that generates the optimal ROC curve of a given score variable directly from its SVT ROC curve, without requiring knowledge of the underlying conditional PDF's. In the case where the PDF's are known, the procedure provides (1) an alternative way of constructing the optimal ROC curve without directly applying LRT's -- apply SVT's then use the procedure to construct the optimal curve -- and (2) an alternative way of finding the score variable values for which the likelihood ratio is greater than a threshold $\eta \geq 0$, without requiring an analytical form for the ratio of the PDF's -- apply SVT's then use the derivative of the SVT ROC curve to read off the values.

\section{Acknowledgements}

The authors would like to thank Professor Meir Feder, Professor George Verghese, Dr. Jim Ward, and Professor Yonina Eldar for their insightful suggestions and discussions about this work.

\bibliographystyle{IEEEbib}
\bibliography{refs}

\end{document}